\magnification=\magstep1
\vbadness=10000
\hbadness=10000
\tolerance=10000

\def\C{{\bf C}}   
\def\L{{\cal L}}  
\def\Q{{\bf Q}}   
\def\R{{\bf R}}   
\def\V{{\cal V}}  
\def\Z{{\bf Z}}   

\proclaim The Gross-Kohnen-Zagier theorem in higher dimensions.
\hfill 1 Oct and 31 Dec 1997, 4 March 2000.

Duke Math. J. 97 (1999), no. 2, 219--233.

Richard E. Borcherds, 
\footnote{$^*$}{Supported by a Royal Society
professorship.}

D.P.M.M.S.,
16 Mill Lane, 
Cambridge, 
CB2 1SB,
England.

e-mail: reb@dpmms.cam.ac.uk

www home page  www.dpmms.cam.ac.uk/\~{}reb

\bigskip 

\proclaim
1.~Introduction.

The Gross-Kohnen-Zagier theorem [G-K-Z] says roughly that the Heegner
divisors of a modular elliptic curve are given by coefficients of a
vector valued modular form of weight $3/2$. We will give another proof
of this (see theorem 4.5 and example 5.1),
which extends to some more general quotients of hermitian
symmetric spaces of dimensions $b^-$ and shows that formal power
series whose coefficients are higher dimensional generalizations of
Heegner divisors are vector valued modular forms of weight $1+b^-/2$.

The main idea of the proof of theorem 4.5 is easy to state. One of the
main results of [B] is a correspondence from modular forms of weight
$1-b^-/2$ with singularities to automorphic forms with known zeros and
poles, which give relations between Heegner divisors and their higher
dimensional generalizations.  On the other hand Serre duality for
modular forms says that the only obstructions to finding modular forms
of weight $1-b^-/2$ with given singularities are given by modular
forms of weight $1+b^-/2$. In other words the only obstructions to
finding relations between Heegner divisors are given by certain
modular forms of weight $1+b^-/2$. It is a formal consequence of this
that the Heegner divisors themselves are the coefficients of a modular
form of weight $1+b^-/2$.  The idea of using the results of [B] to
prove relations between Heegner divisors was suggested to me by
R. L. Taylor.

Most of the more interesting special cases of theorem 4.5 in low
dimensions are already known, though it does at least simplify and
unify several previous proofs of known results.  For modular curves the
theorem is more or less the same as the main result (Theorem C) of
[G-K-Z] stating that Heegner divisors on modular curves are given by
coefficients of a Jacobi form of weight 2. (See example 5.1.) The main
difference is that we prove the result for all Heegner divisors while
the authors of [G-K-Z] restrict to the case of Heegner divisors of
discriminant coprime to the level for simplicity, though their method
could probably be extended to cover all Heegner divisors.  The only
reason this has not been done before (as far as I know) seems to be
that it would take a lot of extra work for a rather small improvement
to the result. Hayashi [H] has extended the results of [G-K-Z] to some
of the other discriminants.

There is a similar result for CM points on Shimura curves (example 5.3).
The abelian varieties in the Jacobians of Shimura curves
are also in the Jacobians of modular curves, so it is likely
that the result for Shimura curves 
is a formal consequence of the result of [G-K-Z] for modular
curves.

In two dimensions theorem 4.5 shows that Hirzebruch-Zagier cycles on a
Hilbert modular surfaces are coefficients of a modular form of weight
$1+b^-/2=2$ (example 5.5).
This is the main result of [H-Z].  There are several other higher
dimensional generalizations of Hirzebruch and Zagier's result about
intersection products of cycles; see for example [K].

In 3 dimensions we recover the result of van de Geer [G] 
that Humbert surfaces on a
Siegel modular variety of dimension 3 are the coefficients of a
modular form of weight $1+b^-/2=5/2$.
There is also a 4 dimensional example due to  Hermann [He],
who showed that the degrees of some 3 dimensional 
modular varieties for the group $Sp_4(\R)$
in a modular variety for the group $U(2,2)$ 
were the coefficients of a modular form of weight $1+b^-/2=3$.

Most previous proofs of the results above prove relations between
divisors by considering an inner product on a space of divisors,
either the N\'eron-Tate inner product in the case of modular curves,
or the intersection number in the case of Hilbert modular
surfaces. The proof in this paper says nothing about inner products,
but instead proves relations between divisors by explicitly
constructing automorphic forms with known poles and
zeros. Unfortunately this means that there is no obvious way to use
the methods of this paper to prove the Gross-Zagier result relating
Heegner divisors to the vanishing of the derivative of an $L$ series.

I would like to thank A. Agboola, J. H. Bruinier, C. Castano-Bernard,
E. Freitag, B. Gross, S. S. Kudla, J. Nekov\'a\v{r},
N. I. Shepherd-Barron, R. L. Taylor, J. Wahl, and D. B. Zagier for
their help. I would particularly like to thank the referee who put a
lot of effort into improving the paper.

%

\proclaim 2.~Modular forms. 

In this section we summarize some standard results about modular forms
and set up notation for the rest of the paper.

Recall that the group $SL_2(\R)$ has a double cover $Mp_2(\R)$
called the metaplectic group whose
elements can be written in the form
$$\left({ab\choose cd},\pm\sqrt{c\tau +d}\right)
$$
where ${ab\choose cd}\in SL_2(\R)$
and $\sqrt{c\tau+d}$ is considered as a holomorphic function
of $\tau$ in the upper half plane whose square is $c\tau+d$.
The multiplication is defined so that the usual formulas for
the transformation of modular forms work for half integer weights,
which means that
$$
(A,f(\cdot))(B,g(\cdot))= (AB,f(B(\cdot))g(\cdot))
$$
for $A,B\in SL_2(\R)$ and $f,g$ suitable functions on $H$.
The group $Mp_2(\Z)$ is the discrete subgroup of $Mp_2(\R)$ 
of elements of the form
$\left({ab\choose cd},\pm\sqrt{c\tau +d}\right)
$ with ${ab\choose cd}\in SL_2(\Z)$. 

Suppose that $\Gamma$ is a subgroup of $Mp_2(\R)$ commensurable 
with $Mp_2(\Z)$, and suppose that $\rho$ is a  representation of $ \Gamma$
on a finite dimensional complex vector space $V_\rho$
which factors through a finite quotient of $\Gamma$
such that $\rho=\sigma_k$ on $\Gamma\cap K$. 
Choose  $k\in {1\over 2}\Z$.
We define a modular form of weight $k$ and type
$\rho$ 
to be a holomorphic function $f$ on the upper half plane $H$
with values in the vector space $V_\rho$ such that
$$f\left({a\tau+b\over c\tau+d}\right)
= \sqrt{c\tau+d}^{2k}
\rho_M\left(\left({ab\choose cd},\sqrt{c\tau+d} \right)\right)
f(\tau)
$$
for  elements $\left({ab\choose cd},\sqrt{c\tau +d}\right)$
of $\Gamma$. (We allow singularities at cusps.)

A modular form has a Fourier expansion at the cusp at infinity as follows.
The Fourier coefficients $c_{n,\gamma}\in \C$ of $f$ are defined by
$$f(\tau) = \sum_{n\in \Q}\sum_{\gamma} c_{n,\gamma}q^ne_\gamma$$
where $q^n$ means $e^{2\pi i n\tau}$ and where the sum runs over a
basis $e_\gamma$ of $V_\rho$ consisting of eigenvectors of $T$. Note
that $n$ is not necessarily integral; more precisely, $c_{n,\gamma}$
is nonzero only if $n\equiv \lambda_\gamma\bmod 1$, where the
eigenvalue of $T$ on $e_\gamma$ is $e^{2\pi i \lambda_\gamma}$.  We
say that $f$ is meromorphic at the cusp $i\infty$ if $c_{n,\gamma}=0$ for
$n<<0$, and we say $f$ is meromorphic at the cusp $a/c$ if
$f((a\tau+b)/(c\tau+d))$ is meromorphic at $i\infty$ for 
${ab\choose cd}\in SL_2(\Z)$. We say that $f$ is holomorphic 
at cusps if the coefficients of the Fourier expansions
at all cusps vanish for $n<0$. We will write 
$ModForm(\Gamma,k,\rho)$ for the space of modular forms of weight 
$k$ and representation $\rho$ for $\Gamma$ which are meromorphic 
at cusps, and 
$HolModForm(\Gamma,k,\rho)$ for the subspace of modular 
forms which are holomorphic
at all cusps.

A particularly important example $\rho_M$ of a representation $\rho$ as above
can be constructed as follows. 
We let $M$ be a nonsingular even lattice of signature $(b^+,b^-)$, with dual
$M'$.  The quotient $M'/M$ is a finite group whose order
is the absolute value of the discriminant of the lattice $M$. 
The mod 1 reduction of $(\lambda,\lambda)/2$ is a
$\Q/\Z$-valued quadratic form on $M'/M$, whose associated
$\Q/\Z$-valued bilinear form is the mod 1 reduction of the bilinear
form on $M'$.  
We let the elements $e_\gamma$ for $\gamma\in M'/M$ be the standard
basis of the group ring $\C[M'/M]$, so that $e_\gamma
e_\delta=e_{\gamma+\delta}$.  
The Grassmannian $G(M)$ of $M$ is defined to be the space of all 
$b^+$-dimensional positive definite subspaces of $M\otimes \R$. 
It is a symmetric space acted on by the orthogonal group
$O_M(\R)$, and if $b^+=2$ it is a hermitian symmetric space.

Recall that there is a unitary
representation $\rho_M$ of the double cover $Mp_2(\Z)$ of $SL_2(\Z)$
on $V_{\rho_M}=\C[M'/M]$ defined by
$$\rho_M(T)(e_\gamma) = e^{2\pi i(\gamma,\gamma)/2}e_\gamma
$$
$$\rho_M(S)(e_\gamma) = {\sqrt{i}^{b^--b^+}\over \sqrt{|M'/M|}}
\sum_{\delta\in M'/M} e^{-2\pi i(\gamma,\delta)}e_\delta
$$
where $T=({11\choose 01},1)$ and $S=({0-1\choose
1\phantom{-}0},\sqrt\tau)$ are the standard generators of $Mp_2(\Z)$,
with $S^2=(ST)^3=Z$, $Z=({-1\phantom{-}0\choose \phantom{-}0-1},i)$,
$\rho_M(Z)(e_\gamma)=i^{b^--b^+}e_{-\gamma}$, $Z^4=1$.  
The representation $\rho_M$
factors through the double cover $Mp_2(\Z/N\Z)$ of the 
finite group $SL_2(\Z/N\Z)$, where $N$ is the smallest integer such that
$N(\gamma,\delta)$ and $N\gamma^2/2$ are integers for all
$\gamma,\delta\in M'$. In particular the representation $\rho_M$
factors through a finite quotient of $Mp_2(\Z)$. 
Also note that there is only one cusp in this case, so that
a modular form is holomorphic if and only if it is holomorphic at $i\infty$. 

\proclaim 3.~Serre duality. 

In this section we show that the space of obstructions 
to finding modular forms with given singularities at cusps 
is dual to a space of holomorphic modular forms. We prove this by identifying
both spaces with cohomology groups of line bundles over modular curves; 
the result we want then follows immediately from  Serre duality.

The space $HolModForm(\Gamma,k,\rho)$ 
can be identified with the space of holomorphic 
sections of the following vector bundle $\V_{k,\rho}$
on $\Gamma\backslash H$, and similarly  $ModForm(\Gamma,k,\rho)$
can be identified with the space of sections which are holomorphic 
outside the cusps. 
We first do the special case when $\rho$ is the trivial representation, 
when $\V_{k,\rho}$ becomes a line bundle $\L_k$. 
Let $K$ be the inverse image of $SO(2)$ in $Mp_2(\R)$, and,
for $k\in {1\over 2}\Z$, define the
character $\sigma_k$ of $K$ by
$$\sigma_k(\theta) = (\pm \sqrt{ci+d})^{2k}$$
for $\theta=({ab\choose cd},\pm\sqrt{c\tau+d})$. 
Note that
$$H \simeq SL_2(\R)/SO(2) \simeq Mp_2(\R)/K,$$
so that there is a homogeneous holomorphic line bundle $\L_k$ on
$H$ defined by
$$\L_k = \bigg(Mp_2(\R)\times \C\bigg)/K,$$
where $\theta\in K$ acts by $(g,z)\theta = (g\theta,\sigma_k(\theta)^{-1}z)$.
For a representation $\rho$ of $\Gamma$ on $V_\rho$
such that $\rho=\sigma_k^{-1}$ on $\Gamma\cap K$, define
a holomorphic vector bundle $\V_{k,\rho}$ over $\Gamma\backslash H$ by
$$\V_{k,\rho} = \Gamma\backslash\bigg( Mp_2(\R)\times V_\rho\bigg)/K,$$
where $K$ acts as before, and $\gamma\in \Gamma$ acts by
$\gamma(g,v) = (\gamma g,\rho(\gamma)v)$.

If $\kappa$ is a cusp of $\Gamma$, let $q_\kappa$ be a uniformizing parameter
at $\kappa$ on $\Gamma\backslash H$. For a representation $\rho$ on
$V_\rho$,
let $V_{\rho}^*$ denote the dual.
Let
$$PowSer_\kappa(\Gamma,\rho)=\C[[q_\kappa]]\otimes V_\rho$$
be the space of formal power series in $q_\kappa$
with coefficients  in $V_\rho$, 
let
$$Laur_\kappa(\Gamma,\rho)=\C[[q_\kappa]][q_\kappa^{-1}]\otimes V_\rho$$
be the space of formal Laurent series in $q_\kappa$
with coefficients  in $V_\rho$, 
and let
$$
Sing_\kappa(\Gamma,\rho)  
= {Laur_\kappa(\Gamma,\rho)\over q_\kappa PowSer_\kappa(\Gamma,\rho)}
$$
be the space of possible 
singularities and constant terms 
of $V_\rho$ valued Laurent series at $\kappa$. 
The two
spaces $PowSer_\kappa(\Gamma,\rho^*)$ 
and $Sing_\kappa(\Gamma,\rho)$ are paired into $\C$ by taking the residue
$$\langle f,\phi\rangle = 
\hbox{Res}\big( f \phi\, q_\kappa^{-1} \,dq_\kappa \big),$$
for $f\in PowSer_\kappa(\Gamma,\rho^*)$ and $\phi\in Sing_\kappa(\Gamma,\rho)$.
Here the product of $f$ and $\phi$ 
is defined using the pairing of $V_\rho$ and $V_\rho^*$.

Then the spaces
$$Sing(\Gamma,\rho) = \oplus_{\kappa} Sing_\kappa(\Gamma,\rho)$$
and
$$PowSer(\Gamma,\rho^*) = \oplus_{\kappa} PowSer_\kappa(\Gamma,\rho^*),$$
where $\kappa$ runs over the $\Gamma$-inequivalent cusps, are paired by the
sum of the  local pairings at the cusps.

There are maps
$$\lambda:HolModForm(\Gamma,k,\rho^*) \longrightarrow PowSer(\Gamma,\rho^*)$$
and
$$\lambda:ModForm(\Gamma,2-k,\rho) \longrightarrow Sing(\Gamma,\rho),$$
defined in the obvious way
by taking the Fourier expansions of their nonpositive part at the
various cusps.

We define the space $Obstruct(\Gamma,k,\rho)$
of obstructions to finding a modular form of type
$\rho$ and weight $k$ which is holomorphic on $H$ and has given
meromorphic singularities and constant terms at the cusps to be the
space
$$Obstruct(\Gamma,k,\rho) =
{Sing(\Gamma,\rho)
\over  \lambda(ModForm(\Gamma,k,\rho))}.
$$

\proclaim Theorem 3.1. 
Suppose that $k\in {1\over 2}\Z$, $\Gamma$ is a subgroup of $Mp_2(\R)$ which
is commensurable with $Mp_2(\Z)$ and $\rho$ is a finite dimensional
complex representation of $\Gamma$ factoring through a finite quotient of
$\Gamma$ such that $\rho=\sigma_k$ on $\Gamma\cap K$.  
Then the space of obstructions $Obstruct(\Gamma,2-k,\rho)$
is finite dimensional and 
dual to the space $HolModForm(\Gamma,k,\rho^*)$. 
The pairing between them is induced by the above pairing
between $Sing(\Gamma,\rho)$ and $PowSer(\Gamma,\rho^*)$. 
In other words, 
$$\lambda\bigg(ModForm(\Gamma,2-k,\rho)\bigg) =
\lambda\bigg(HolModForm(\Gamma,k,\rho^*)\bigg)^\perp,$$
and also, since the pairing is nondegenerate, 
$$\lambda\bigg(HolModForm(\Gamma,k,\rho^*)\bigg) =
\lambda\bigg(ModForm(\Gamma,2-k,\rho)\bigg)^\perp.$$

Proof.  Suppose first that the group $\Gamma$
acts freely on the upper half plane $H$, and the representation $\rho$
is one dimensional and trivial.  For
any  $k\in {1\over 2}\Z$ we let $\L_k$ be the line bundle defined above
whose sections are holomorphic  modular forms of weight $k$.
We let $\L_{cusp}$ be the line
bundle corresponding to the divisor which is the union of all
cusps of the compactification of $\Gamma\backslash H$.
The canonical line bundle  is isomorphic to $\L_2\otimes
\L_{cusp}^*$ because holomorphic 1-forms are essentially the same as
cusp forms of weight 2.  By Serre duality we see that the space
$H^0(\L_k)$ of holomorphic forms of weight $k$ is dual to
$H^1(\L_2\otimes \L_{cusp}^*\otimes \L_k^*)= H^1(\L_{2-k}\otimes
\L_{cusp}^*)$.  
If $\L$ is any line bundle over  a compact Riemann surface then
the cohomology group $H^1(\L)$ can be identified with the space of
obstructions to finding a meromorphic section of $\L$ with given
singularities at some fixed nonempty finite set of points and
holomorphic elsewhere.  Hence $H^1(\L_{2-k}\otimes \L_{cusp}^*)$ is
the space of obstructions to finding a meromorphic section of
$\L_{2-k}\otimes \L_{cusp}^*$ with given singularities at cusps and
holomorphic elsewhere (since there is at least one cusp).  This in
turn is the space of obstructions $Obstruct(\Gamma, 2-k, \C)$
to finding a meromorphic section of
$\L_{2-k}$ with given singularities and constant terms at cusps and
holomorphic elsewhere.  The pairing in Serre duality is the one above 
given by
taking the sum of residues of the product.
This proves theorem 3.1 in the case that
$\Gamma$ acts trivially on $\rho$ and fixed point freely on $H$.

In the general case choose a finite index subgroup $\Gamma_0$ of
$\Gamma$ such that $\Gamma_0$ acts trivially on $\rho$ and acts fixed
point freely on $H$. Then the cohomology 
groups $H^1$ and $H^0$ for $\Gamma$ are
just the fixed points under $\Gamma/\Gamma_0$ of the corresponding
cohomology groups for $\Gamma_0$. The theorem for $\Gamma$ now follows from the
fact that if we have a perfect duality between finite dimensional complex
vector spaces that is invariant under some finite group
$\Gamma/\Gamma_0$ acting on these vector spaces, then we also get a
perfect duality between the fixed points of $\Gamma/\Gamma_0$
(as follows from the fact that finite dimensional complex representations
of finite groups are completely reducible).  This
proves theorem 3.1.

\proclaim 4.~Heegner divisors as Fourier coefficients. 

In this section we prove the main theorem 4.5, saying roughly that
Heegner divisors are the Fourier coefficients of a modular form.
We do this by using the modular forms with singularities constructed
in section 3 to find a large number of relations between Heegner divisors.

We let $M$ be an even lattice of signature $(2,b^-)$. 
Suppose that $\Gamma$ is a discrete group acting on the Grassmannian
$G(M)$ of $M$. We define a divisor on $X_\Gamma=\Gamma\backslash G(M)$ to
be a locally finite $\Gamma$ invariant divisor on $G(M)$ whose support
is a locally finite union of a finite number of $\Gamma$-orbits of
irreducible codimension 1 subvarieties of $G(M)$. 
This definition is
a sort of crude substitute for the definition of a divisor on the
``orbifold'' (or algebraic stack) $X_\Gamma$ but is adequate
for the purposes of this paper.  Note that this is not quite the same as a
divisor on the complex analytic space $X_\Gamma$; for example,
if $G(M)$ is the upper half plane $H$ and $\Gamma$ is $SL_2(\Z)$ then
the image of the point $i$ represents a divisor in $\Gamma\backslash
H$, which is twice a divisor in the orbifold $\Gamma\backslash H$ but
not in the complex manifold $\Gamma\backslash H$.

Recall that for any negative norm vector $v\in M\otimes \R$ there is a
divisor $v^\perp$ of $G(M)$, equal to the points of the Grassmannian
represented by 2-planes orthogonal to $v$.  If $n$ is any negative
rational number and $\gamma\in M'/M$ then we define the Heegner
divisor $y_{n,\gamma}$ to be the sum of the divisors of all norm $2n$
vectors of $M+\gamma$.  Note that $y_{n,\gamma}=y_{n,-\gamma}$ because
if $v\in M+\gamma$ then $-v\in M-\gamma$.  We define the group
$Heeg(X_\Gamma)$ of Heegner divisors to be the direct sum of a copy of
$\Z$ generated by a symbol $y_{0,0}$, and the subgroup of the group of
divisors generated by the Heegner divisors $y_{n,\gamma}$.  If $n>0$
or $n=0,\gamma\ne0$ then we define $y_{n,\gamma}$ to be 0.  We define
a Heegner divisor to be principal if it is of the form $c_{0,0}y_{0,0}
+ D$, where $D$ is the divisor of a meromorphic automorphic form of
weight $c_{0,0}/2$ for some integer $c_{0,0}$ and some unitary
character of finite order of the subgroup ${\rm Aut}(M)$ fixing all
elements of $M'/M$. Here the weight of an automorphic form is the
weight used in theorem 13.3 of [B]. We write $PrinHeeg(X_\Gamma)$ for
the subgroup of principal Heegner divisors, and $HeegCl(X_\Gamma)$ for
the group $Heeg(X_\Gamma)/PrinHeeg(X_\Gamma)$ of Heegner divisor
classes.  (The published version of this paper omitted the condition
that the character have finite order. J. Bruinier pointed out to me
that there are sometimes ``too many'' automorphic forms with
characters of infinite order (see [F]), so that if infinite order
characters are allowed the group of Heegner divisor classes would
sometimes collapse.)

There is a surjective linear map
$$
\xi:Sing(Mp_2(\Z),\rho_M) \longrightarrow Heeg(X_\Gamma)\otimes_\Z\C
$$
taking $q^ne_\gamma$ to $y_{n,\gamma}$. 
For a subring $F$ of $\C$ let
$Sing(Mp_2(\Z),\rho_M)_F$ be the $F$-submodule of $Sing(Mp_2(\Z),\rho_M)$
for which the coefficients of $q^ne_\gamma$ for $n\le 0$ are in $F$, and let
$$ModForm(Mp_2(\Z),1-b^-/2,\rho_M)_\Z 
\subseteq ModForm(Mp_2(\Z),1-b^-/2,\rho_M)$$
be the $\Z$-submodule
whose image under $\lambda$ lies in $Sing(Mp_2(\Z),\rho_M)_\Z$.

\proclaim Theorem 4.1.   
Suppose $M$ is an even lattice of signature $(2,b^-)$ and $f$ is a
modular form of weight $1-b^-/2$ and representation $\rho_M$ which is
holomorphic on $H$ and meromorphic at cusps and whose coefficients
$c_{n,\gamma}$ are integers for $n\le 0$.  Then $\sum_{n,\gamma}
c_{n,\gamma}y_{n,\gamma}$ is a principal Heegner divisor.
In other words, 
$$\xi\bigg(\lambda\big(ModForm(Mp_2(\Z),1-b^-/2,\rho_M)_\Z\big)\bigg) \subseteq
PrinHeeg(X_\Gamma).$$

Proof. 
Theorem 13.3 of [B] implies that if $M$ and $f$ satisfy the conditions 
of theorem 4.1 
then there is a meromorphic  function $\Psi$ on $G(M)$
with the following properties.
\item{1.} $\Psi$ is an automorphic form of weight $c_{0,0}/2$
for some unitary character of the subgroup of ${\rm Aut}(M)$
fixing all elements of $M'/M$.
\item{2.} The only zeros or poles  of $\Psi$ lie on the
divisors $\lambda^\perp$
for $\lambda\in M$, $\lambda^2<0$ and are zeros of order
$$\sum_{0<x\in \R\atop x\lambda\in M'}c_{x^2\lambda^2/2,x\lambda}
$$
(or poles if this number is negative).

Theorem 4.1 is just a restatement of this result 
in terms of the divisors $y_{n,\gamma}$, 
provided we show that the characters 
of the automorphic forms have finite order. This can be shown as follows. 
For $O_{2,n}(\R)$ with $n >2$ this follows because these Lie groups have no
almost simple factors of real rank 1, and if $G$ is a lattice in a
connected Lie group with no simple factors of rank 1 then the
abelianization of $G$ is finite.  (See [M page 333,  proposition 6.19].)
So any character of $G$ has finite order.

For the cases $n=1$ and $n=2$ we use the embedding trick ([B98, lemma
8.1]) to see that if $f$ is an infinite product of $O_{2,n}(\R)$ then
$f$ is the restriction of an infinite product $g$ of 
$O_{2, 24+n}(\R)$. The infinite product $g$ is not necessarily single
valued; however, a look at the proof of lemma 8.1 shows that if $f$ is
constructed from a vector valued modular form with integral
coefficients, then $g^{24}$ has zeros and poles of integral order and
is therefore a meromorphic automorphic form for some unitary
character. By the previous paragraph this character has finite order, 
and therefore so does the character of $f$. 
This proves theorem 4.1.

We remark that the proof of theorem 13.3 of [B] is quite long but
most of the proof is the calculation of the Fourier expansion of
$\Psi$ (or rather of its logarithm), which we do not use in this
paper.  The proof of the parts of theorem 13.3 that we use here is
much shorter and is mainly in section 6 of [B].

\proclaim Lemma 4.2.  There is a number field $F$ of finite
degree over $\Q$ such that the finite dimensional space
$HolModForm(Mp_2(\Z),1+b^-/2,\rho_M^*)$
has a basis whose Fourier coefficients all lie in $F$, i.e.,
such that $\lambda(f)\in PowSer(Mp_2(\Z),\rho_M^*)_F$. 

Proof. It follows from [S, section 3.5] that the space of
modular forms of level $N$ has a base of forms whose Fourier 
expansions at all cusps have coefficients in some algebraic number field
of finite degree. (Shimura covers the case of integral weight at least 2,
but the other cases can easily be reduced to this by multiplying
forms by a power of $\eta(\tau)$; note that if we are given a 
basis for the space of modular forms of level $N$ all of whose Fourier
expansions have coefficients in $F$ then we can find a similar basis
for the space of forms which vanish to given orders at various cusps.)
Each of the $|M'/M|$ components of $HolModForm(Mp_2(\Z),1+b^-/2,\rho_M^*)$
is a modular form of level $N$ for some $N$, and the representation
$\rho_M$ is obviously
defined over an algebraic number field of finite degree. 
This implies that the space $HolModForm(Mp_2(\Z),1+b^-/2,\rho_M^*)$
has a basis of elements all of whose coefficients lie in 
some finite algebraic number field, because this space is 
the $Mp_2(\Z/N\Z)$-invariant subspace of the space of modular
forms of level $N$ with coefficients in $\C[M'/M]^*$. 
This proves lemma 4.2.

\proclaim Lemma 4.3. 
Let $Gal(\bar{\Q}/\Q)\cdot\lambda(HolModForm(Mp_2(\Z),1+b^-/2,\rho_M^*))$ be
the
space of $Gal(\bar{\Q}/\Q)$ conjugates of the $q$-expansions of elements of
$HolModForm(Mp_2(\Z),1+b^-/2,\rho_M^*)$. (This is well defined
and finite dimensional by lemma 4.2.) Then
$$\eqalign{
&\lambda\big(ModForm(Mp_2(\Z),1-b^-/2,\rho_M)_\Z\big)\otimes\C\cr
=&
\bigg(Gal(\bar{\Q}/\Q)\cdot
\lambda\big(HolModForm(Mp_2(\Z),1+b^-/2,\rho_M^*)\big)\bigg)^\perp.\cr
}$$
Moreover, this space has finite index in
$Sing(Mp_2(\Z),\rho_M)$.

Proof. It is obvious using theorem 3.1 
that the first space is contained in the second. 
To show that the second space is contained in the first, 
choose a basis for $HolModForm(Mp_2(\Z),1+b^-/2,\rho_M^*)$ consisting
of forms with coefficients in some algebraic number field $F$ of finite degree,
which we can do by lemma 4.2. But then the image of this space under
$Gal(\bar \Q/\Q)$ is spanned by a finite number of functions all of whose
Fourier coefficients are rational. So the orthogonal complement
is defined by rational linear relations and therefore has
a basis consisting of modular forms whose image under $\lambda$
has rational coefficients. By multiplying by a common denominator we can 
assume the coefficients are integral. By applying 3.1 again 
this proves that the second space is
contained in the first, so both spaces are equal. 

The fact that the space has finite index in $Sing(Mp_2(\Z),\rho_M)$
follows because it is the orthogonal complement of the sum of a finite
number of $Gal(\bar \Q/\Q)$ conjugates of a finite dimensional
space. This
proves lemma 4.3.

\proclaim Lemma 4.4. The complex vector space $HeegCl(X_\Gamma)\otimes\C$
generated by the Heegner divisor classes
is finite dimensional.

{Proof.}
By Theorem~4.1, the surjective map
$$Sing(Mp_2(\Z),\rho_M) 
\longrightarrow Heeg(X_\Gamma)\otimes_\Z\C\longrightarrow
HeegCl(X_\Gamma)\otimes\C$$
factors through the  space
$$Sing(Mp_2(\Z),\rho_M)/
\lambda(ModForm(Mp_2(\Z),1-b^-/2,\rho_M)_\Z)\otimes\C.$$
which is finite dimensional by lemma 4.3. 
This proves lemma 4.4.

\proclaim Theorem 4.5. If the Heegner divisors $y_{n,\gamma}$
are considered as elements of the Heegner divisor class group
$HeegCl(X_\Gamma)\otimes \C$
then
$$\sum_{n\in \Q} \sum_{\gamma\in M'/M} y_{-n,\gamma} q^ne_\gamma$$ is
a modular form, and more precisely it lies in the space 
$$
\bigg(HeegCl(X_\Gamma)\otimes\C\bigg)\otimes_\C
\bigg(Gal(\bar{\Q}/\Q)\cdot
\lambda\big(HolModForm(Mp_2(\Z),1+b^-/2,\rho_M^*)\big)\bigg).
$$

Proof. 
By lemma 4.4 we know that 
$$\sum_{n\in \Q}\sum_{\gamma\in M'/M} y_{-n,\gamma} q^n e_\gamma \in
\bigg(HeegCl(X_\Gamma)_\C\bigg)\otimes_\C PowSer(Mp_2(\Z),\rho_M^*).$$
By Theorem~4.1, the pairing 
$$\sum c_{n,\gamma}y_{n,\gamma} = \xi(\lambda(\sum c_{n,\gamma} q^ne_\gamma))$$
of this series with any
element $\sum c_{n,\gamma}q^ne_\gamma$ of
$$\lambda\big(ModForm(Mp_2(\Z),1-b^-/2,\rho_M)_\Z\big)$$
is zero. Theorem 4.5 therefore follows from lemma 4.3 above,

 In the few examples I have checked the space of
modular forms of weight $1+b^-/2$ and type $\rho_M^*$
has a basis of forms with rational coefficients, so that the Heegner 
divisors are the coefficients of a modular form of type $\rho_M^*$.  
I do not know whether or not such a basis always exists. 

 Modular forms of type $\rho_M^*$ and weight $1+b^-/2$ can be
identified with certain Jacobi forms in several variables of weight
$1+b^-$ as in [E-Z theorem 5.1], so theorem 4.5 says that the Heegner
divisors are coefficients of a Jacobi form of weight $1+b^-$.

We can often work out the dimensions of the spaces of vector valued
modular forms using the Riemann-Roch theorem. For example,
if $\rho$ is a $d$-dimensional representation of $Mp_2(\Z)$ on which $Z$ acts
as $e^{-\pi i k}$ for some $k\ge 2$ with $k\in {1\over 2}\Z$ 
then the dimension of the space
of holomorphic modular forms of type $\rho$ and weight $k$ is equal to
$$d+dk/12 -\alpha(e^{\pi i k/2}S) -\alpha((e^{\pi i k/3}ST)^{-1})-\alpha(T)$$
where 
$\alpha(X)$ is the is the sum of the numbers $\beta_j$, $1\le j\le d$,
where the eigenvalues of $X$ 
are $e^{2\pi i \beta_j}$ and $0\le \beta_j <1$. 
This formula cannot be applied directly to $\rho_M^*$
as the condition on $Z$ is not satisfied, but we can apply it to 
the subspace of $\rho_M^*$ on which $Z$ acts as $e^{-\pi i k}$; 
for $k=1+b^-/2$
this is the subspace spanned by the elements $e^*_\gamma+e^*_{-\gamma}$.
The equality $y_{n,\gamma}=y_{n,-\gamma}$ implies that the modular
form in theorem 4.5 lies in the space 
on which $Z$ acts as $e^{-\pi i k}$.

\proclaim 5.~Examples.

In this section we give some examples to illustrate
theorem 4.5. In most cases we will describe the lattice $M$
and the Heegner divisors. 

{\bf Example 5.1} We work out theorem 4.5 in the case of modular curves;
see [Z84] and [G-K-Z] for more about this case. 
We fix $N$ to be any positive integer (called the level). 

We let $M$ be the 3 dimensional even lattice of all symmetric 
matrices $v=\pmatrix{C/N&-B/2N\cr -B/2N&A/N}$ with $A/N, B/2N, C$ integers,
with the norm $(v,v)$ defined to be $-2N\det(v) = (B^2-4AC)/2N$. 
The dual lattice is the set of matrices as above with $A/N,B,C\in \Z$, 
and   $M'/M$ can be identified with $\Z/2N\Z$
by mapping a matrix of $M'$ to the value of $B\in \Z/2N\Z$. 
The lattice $M$ splits as the direct sum of
the 2 dimensional hyperbolic unimodular even lattice $II_{1,1}$ and
a lattice generated by an element of norm $2N$. 
The group $\Gamma_0(N)=\{{ab\choose cd}\in SL_2(\Z)|c\equiv 0\bmod N\}$
acts on the lattice $M$ by $v\mapsto XvX^t$ for
$X\in \Gamma_0(N)$, and under this action it fixes all elements of $M'/M$. 

We identify the upper half plane with points in the 
Grassmannian $G(M)$
by mapping
$\tau\in H$ to the 2 dimensional positive definite space spanned by
the real and imaginary parts of the norm 0 vector 
$\pmatrix{\tau^2&\tau\cr\tau&1\cr}$.
For each  $n\in \Q$ and $\gamma\in M'/M=\Z/2N\Z$  the Heegner
divisor $y_{n,\gamma}$ is the union of the points orthogonal to
norm $2n$ vectors of $M+\gamma$.
In terms of points on $H$ this Heegner
divisor consists of all points $\tau\in H$ such that
$$A\tau^2+B\tau+C=0$$ for some integers $A$, $B$, $C$ (not necessarily
coprime) with $N|A$, $B\equiv \gamma\bmod 2N$, $B^2-4AC=4Nn$.  (Warning:
we omit the condition $A>0$ which is sometimes included in the
definition of Heegner divisors. The points with $A>0$ are exchanged
with the points with $A<0$ under the Fricke involution $\tau\mapsto-1/N\tau$,
so the
difference is not important as long as we quotient out by the Fricke
involution.)  The Heegner divisor $y_{n,\gamma}$ is invariant under
the group $\Gamma_0(N)$ and the Fricke involution,
because the corresponding set of vectors is invariant under
$\Gamma_0(N)$ and changes sign under the Fricke involution.

The element $y_{0,0} $ is 0 (modulo torsion) as follows easily from 
the existence of the modular form $\Delta$ with no zeros on $H$, and
this is why $y_{0,0}$ does not  appear explicitly in this case.
However $y_{0,0}$ is usually nonzero in other cases; see example 5.4. 
The vanishing of $y_{0,0}$ seems to be related to the failure of the
Koecher boundedness principle. 

To compare theorem 4.5 with theorem C of [G-K-Z] which is stated in
terms of Jacobi forms rather than vector valued modular forms we
recall that according to [E-Z, theorem 5.1] the space $J_{k,N}$ of
Jacobi forms of weight $k$ and index $N$ is isomorphic to the space of
modular forms of weight $k-1/2$ and type $\rho_M^*$. By [E-Z theorem
9.3] the spaces of Jacobi forms $J_{2,N}$ have bases consisting of
forms with rational coefficients, so by the remarks after theorem 4.5
the series
$$ \sum_{n,\gamma}y_{-n,\gamma}q^ne_\gamma$$ is a modular form of
weight $3/2$ and type $\rho_M^*$.  This almost implies the main result
theorem C of [G-K-Z], and (assuming the Gross-Zagier theorem [G-Z]) is
essentially the ``ideal statement'' of theorem C conjectured on page
503 of [G-K-Z]. (But note that there is a small technical difference
between our definitions and those of [G-Z], because we allow
nontrivial unitary characters of finite order in the definition of
principal Heegner divisors. This makes no difference, because some
finite positive power of an automorphic form with character of finite
order has trivial character, 
so our definitions are equivalent to those of [G-Z].)

We have implicitly quotiented out by the divisors of cusps.
We can ask what happens if we leave them in, which essentially 
comes down to looking at the power series whose coefficients
are given by the degrees of Heegner divisors. These degrees are 
given by various sorts of class numbers (e.g. Hurwitz class number in 
the case of level 1) and Zagier showed ([Z], [H-Z chapter 2]) 
that these are the 
coefficients of certain non holomorphic modular forms of weight 3/2. 

{\bf Example 5.2.} We can ask if all relations between Heegner
divisors are given by theorem 4.5. The answer seems to be that usually
they are but there are some exceptions where there are extra relations.
(Remark added 2000: this has mostly been proved by J. Bruinier [Br].)
For example, by the Gross-Zagier theorem [G-Z] this happens whenever
there is a modular elliptic curve of odd rank at least 3, because for
such elliptic curves all Heegner points vanish. For $b^->2$ it might
be possible to prove that there are no further relations between
Heegner points as follows. For each Heegner point we can construct a
real analytic function with a logarithmic singularity along the
Heegner divisor by applying the singular theta correspondence to a
possibly non holomorphic modular form. For any linear combination of
Heegner divisors we get a function with logarithmic singularities
along these divisors, which is the logarithm of an automorphic form
when there is a relations between the Heegner divisors given by
theorem 4.5. When there is no relation it should be possible to prove
that there is no automorphic form giving this relation, as otherwise
the difference of the log of this form and the function we have
constructed would be a nonzero harmonic function vanishing at
infinity, which is not possible. This argument breaks down for $b^-=$
1 or 2 because in attempting to construct non holomorphic modular
forms with given singularities at cusps we run into problems with the
Green's functions having poles in the critical strip. These poles are
presumably related to the derivative of $L$ functions at $s=1$, and it
might be possible to prove a weak version of the Gross-Zagier theorem
along these lines, saying that certain spaces spanned by Heegner
points have rank 0 if and only if the derivative of some $L$ series
vanishes at $s=1$.

{\bf Example 5.3.} We can do the same as in example 5.1 with the
Shimura curves associated with quaternion algebras over the
rationals. Suppose that $K$ is a non split 4 dimensional central
simple algebra over the rationals which is split at infinity and $R$
is some order in $K$.  We let $M$ be the lattice of elements of $R$
orthogonal to 1, with the inner product given by minus that of $R$.
Then $M$ is a lattice in $\R^{2,1}$, and the group of units of $R$ acts
on $M$ by conjugation. The Grassmannian $G(M)$ is isomorphic to the
upper half plane $H$ and the group $\Gamma$ acts on it with the
quotient being a compact Riemann surface (or more precisely an
orbifold) called a Shimura curve.  The points on the Shimura curve
associated to vectors of $M$ are called CM points. So theorem 4.5
implies that certain divisors associated to CM points are coefficients
of modular forms of weight $3/2$.
The case of modular curves for $\Gamma_0(N)$ is really a special case
of the construction above where we take $K$ to be the split central
simple algebra $M_2(\Q)$.

{\bf Example 5.4.} In the case of the Gross-Kohnen-Zagier theorem the
modular form we get is a cusp form and $y_{00}=0$.  We give an example
to show that in higher dimensions the modular forms we get are not
necessarily cusp forms and $y_{00}$ can be nonzero.  We take $M=M'$ to
be the lattice $II_{25,1}$. Then for each positive  integer $n$ we
have a Heegner divisor $y_{-n,0}$ (which is irreducible if and only if
$n$ is square free).  The representation $\rho_M$ is trivial, so
vector valued modular forms of type $\rho_M^*$ are the same as modular
forms of level 1 and in particular we can find a base of forms with
rational coefficients.  By theorem 4.5 the sum $\sum q^ny_{-n,0}$ is 
then a modular form
of weight 14, so must be a multiple of
$E_{14}=1-24\sum_n\sigma_{13}(n)q^n$.  Therefore $y_{-n,0} =
-24\sigma_{13}(n) y_{0,0}$. The Koecher boundedness principle
implies that no positive integral multiple of any $y_{-n,0}$ for any
positive integer $n$ is zero.  This example also shows that $y_{0,0}$
is not necessarily in the integral lattice generated by Heegner
divisors.  The automorphic form giving the relation $y_{-1,0} = -24
y_{0,0}$ is unusual because it is a holomorphic form of singular
weight and is also the denominator function of the fake monster Lie
algebra.

{\bf Example 5.5.} Hirzebruch and Zagier showed in [H-Z] that a power
series with coefficients given by modular curves on a Hilbert modular
surface was a modular form of weight 2. We will show that similar
results can be deduced from theorem 4.5. (It seems likely that these
results imply Hirzebruch and Zagier's result, but I have not checked
this in detail.)

For background about Hilbert modular surfaces see [G] or [H-Z].  We
fix a real quadratic field $K$ of discriminant $D>0$ and ring of
integers $R$. We denote the conjugate of $\lambda \in K$ by
$\lambda'$.  We will do the case of Hilbert modular surfaces of level
1 associated to the trivial ideal class; it should be obvious how to
choose $M$ to extend this to higher levels and other ideal classes.

We let $M$ be the even lattice of matrices of the form $v=\pmatrix{
C&-B\cr -B'&A\cr}$ with $A,C\in \Z$, $B\in R$, with
norm given by $-2\det(v)$. 
The group $SL_2(K)$ acts on
the vector space $M\otimes \Q$ of hermitian matrices
by $v\rightarrow XvX'^t$ for $X\in
SL_2(K)$ and $v\in M\otimes \Q$.  The group $SL_2(R)$ maps $M$
to itself under this action.  The lattice $M$ is the direct sum of the
lattice $II_{1,1}$ and a lattice of determinant $D$.

We identify the product of two copies of the upper half plane $H$ with
the Grassmannian of $M$ by mapping $(\tau_1,\tau_2)\in H^2$ to the
space spanned by the real and imaginary parts of the norm 0 vector
$\pmatrix{\tau_1\tau_2&\tau_1\cr \tau_2&1\cr}$.  This
induces the usual action of $SL_2(K)$ on $H^2$ given by
$${ab\choose cd}((\tau_1,\tau_2))
=\left({a\tau_1+b\over c\tau_1+d},{a'\tau_2+b'\over c'\tau_2+d'}\right).$$

If $v$ is a negative norm vector of $M'$ then we define the modular
curve $T_v$ to be the orthogonal complement of $v$ in the Grassmannian
of $M$. We can describe $T_v$ explicitly as follows.  If $v$ is the
matrix $\pmatrix {C&-B\cr -B'& A\cr}$ then $T_v$ is the set
of points $(\tau_1,\tau_2)\in H^2$ such that
$$A\tau_1\tau_2 + B'\tau_1+B \tau_2 + C = 0.$$ If
$0>n\in \Q$ and $\gamma\in M'/M$ then $y_{n,\gamma}$ is
the union of all the curves $T_v$ for norm $2n$ vectors $v\in
M+\gamma$. Then theorem 4.5 implies that for some choice of $y_{0,0}$
the power series
$$\sum_{n,\gamma} y_{-n,\gamma}q^ne_\gamma$$
is a modular form of weight 2.

In [H-Z, Theorem 1] Hirzebruch and Zagier show that for certain
divisors $T_N$ for $N>0$ the power series $\sum_N T^c_Nq^N$ is a modular
form of weight 2. (Here $T^c_N$ is the homology class 
of the Hilbert modular surface given by the projection
of the homology class of $T_N$ into the orthogonal complement of the subspace
generated by homology cycles of the curves of the cusp resolutions.)
This result is closely related to the one
above, because the divisors $T_N$ are the unions of the divisors
$y_{-N,\gamma}$ for all $\gamma\in M'/M$.  
There are several minor differences between the results
here and the result in [H-Z], as follows. For simplicity Hirzebruch
and Zagier only treat the case of $K$ having discriminant $D$ a prime
$p$ congruent to 1 mod 4, but their methods could probably
be extended
to all positive discriminants.  The lattice used by [H-Z] is $pM'$ of
discriminant $p^3$ rather than $M$ of discriminant $p$. 
The definition of a principal Heegner divisor here is slightly 
different because we allow nontrivial unitary characters
(though it it likely that all these characters have finite order, 
in which case there is no essential difference).  There are
also some other small changes of notation; for example, we use
hermitian matrices rather than skew hermitian matrices to emphasize
the similarity with example 5.1.

This example is closely related to the case of modular curves in 5.1;
in fact in [Z84] Zagier shows how to deduce the results about modular
curves from the result on Hilbert modular surfaces in many special cases.

{\bf Example 5.6.} If we take $M$ to be a 5 dimensional lattice then
the symmetric space of $M$ is isomorphic to the Siegel upper half
plane of genus 2. The divisors on this Siegel upper half plane
associated to vectors of $M$ (or rather their images in the quotient)
are just the so-called Humbert surfaces. Theorem 4.5 implies that power
series with coefficients given by certain Humbert divisors are 
vector valued modular
forms of weight 5/2. A similar result is mentioned in [G, p. 213].

{\bf Example 5.7.} If we take $M$ to be a 6 dimensional lattice 
then the symmetric space of $M$ is isomorphic to 
the Hermitian upper half space of complex dimension 4. The divisors
associated to vectors of $M$ are quotients of the Siegel 
upper half space of complex dimension 3.  Hermann [He] found an
example of a modular form of weight $1+b^-/2=3$ associated to this case. 
Note that the groups $Sp_4(\R)$ and $SU(2,2)$ used in [He]
are locally isomorphic
to the groups $O_{2,3}(\R)$ and $O_{2,4}(\R)$ used in this paper,
and in particular the Hermitian half space of [He] is
isomorphic to the Grassmannian of $\R^{2,4}$.

\proclaim References.

\item{[B]} R. E. Borcherds, 
Automorphic forms with singularities on Grassmannians,
alg-geom/9609022. Accepted by Invent. Math. 
\item{[Br]}
J. H. Bruinier, Borcherds products and Chern classes of
Hirzebruch-Zagier divisors. Invent. Math. 138 (1999), no. 1, 51--83.
Borcherds products on O(2,l) and Chern classes of Heegner divisors,
1999 preprint, available from {\tt
http://www.rzuser.uni-heidelberg.de/\hbox{\~{}}t91/}
\item{[E-Z]}{M. Eichler, D. B. Zagier, 
``The theory of Jacobi forms'', Progress
in mathematics vol. 55,
Birkh\"auser, Boston, Basel, Stuttgart 1985.}
\item{[F]}
J. D. Fay,  Fourier coefficients of the resolvent for a Fuchsian
group. J. Reine Angew. Math. 293/294 (1977), 143--203.
\item{[G]} G. van der Geer, ``Hilbert modular surfaces'', 
Springer Verlag 1987.
\item{[G-K-Z]} B. H. Gross, W. Kohnen, D. B. Zagier, Heegner points and 
derivatives of $L$-series. II. Math. Ann. 278
(1987), no. 1-4, 497--562.
\item{[G-Z]}B. H. Gross,  D. B. Zagier,  Heegner points and derivatives
of $L$-series. Invent. Math. 84
(1986), no. 2, 225--320.
\item{[H]} Y. Hayashi, 
The Rankin's $L$-function and Heegner points for general discriminants. 
Proc. Japan Acad. Ser. A Math. Sci. 71 (1995), no. 2, 30--32. 
\item{[He]} C. F. Hermann, Some modular varieties related to $P^4$. 
In ``Abelian varieties'', edited by Barth, Hulek, and Lange, 
Walter de Gruyter Verlag, Berlin--New York, (1995). 
\item{[H-Z]} F. Hirzebruch, D. B. Zagier, 
Intersection numbers of curves on Hilbert modular surfaces and modular
forms of Nebentypus. Invent. Math. 36 (1976), 57--113.
\item{[K]} S. S. Kudla, Algebraic cycles on Shimura varieties of
orthogonal type, Duke Math. J. Vol 86 No. 1 (1997), 39--78.
\item{[M]} 
G. A. Margulis,  Discrete subgroups of semisimple Lie
groups. Ergebnisse der Mathematik und ihrer Grenzgebiete (3),
17. Springer-Verlag, Berlin, 1991. ISBN: 3-540-12179-X
\item{[S]} G. Shimura, Introduction to the arithmetic theory of automorphic 
functions, Princeton University Press, 1970. 
\item{[Z]} D. B. Zagier, 
Nombres de classes et formes modulaires de poids $3/2$,
C. R. Acad. Sci. Paris S\'er. A-B 281 (1975), no. 21, Ai, A883--A886.
\item{[Z84]} D. B. Zagier,
Modular points, modular curves, modular surfaces and modular forms. 
Workshop Bonn 1984 (Bonn, 1984), 225--248, 
Lecture Notes in Math., 1111, Springer, Berlin-New York, 1985.

\bye